\documentclass{article}

\usepackage {amsmath, amssymb, amsthm}
\usepackage [top=1.5in, right=1in, bottom=1.5in, left=1in] {geometry}
\usepackage {pgfplots}
\pgfplotsset {compat=1.13}
\usepackage {tikz}
\usepackage{todonotes}
\usepackage{cite}
\usepackage{graphicx}
\usepackage {hyperref}
\hypersetup {
    colorlinks = true,
    linktoc = all,
    linkcolor = blue,
    citecolor = red
}

\usepackage {libertine}
\usepackage [T1]{fontenc}

\theoremstyle {plain}
\newtheorem {thm} {Theorem}[section]
\newtheorem {prop} [thm]{Proposition}

\newtheorem {lem}[thm] {Lemma}

\newtheorem {problem} [thm]{Problem}

\theoremstyle {definition}
\newtheorem {rmk}[thm] {Remark}

\newcommand{\R}{\mathbb{R}}

\newcommand{\Z}{\mathbb{Z}}

\newcommand{\trop}{\mathfrak{m}}

\DeclareMathOperator{\sign}{sign}

\title{\vspace{-2.5cm}Discrete Dynamical Systems From Real Valued Mutation}
\author{John Machacek$^1$ and Nicholas Ovenhouse$^2$}
\date{$^1$Hampden Sydney College \\
      $^2$University of Minnesota}

\begin {document}

\maketitle

\begin{abstract}
    We introduce a family of discrete dynamical systems which includes, and generalizes,
    the mutation dynamics of rank two cluster algebras. These systems exhibit behavior associated with integrability,
    namely preservation of a symplectic form, and in the tropical case, the existence of a conserved quantity. 
    We show in certain cases that the orbits are unbounded.
    The tropical dynamics are related to matrix mutation, from the theory of cluster algebras. 
    We are able to show that in certain special cases, the tropical map is periodic.
    We also explain how our dynamics imply the asymptotic sign-coherence observed by Gekhtman and Nakanishi in the $2$-dimensional situation.
\end{abstract}

\tableofcontents

\section {Introduction}

For any two positive real numbers $p$ and $q$, we define two discrete dynamical systems inspired by cluster algebra theory.
The first system maps the first quadrant to itself in a way which generalizes rank two cluster algebra mutation to matrices with real entries.
It is given by iterating the following map
\begin {equation} \label{eq:mu}
    \mu(x,y) = \left( \frac{1+y^q}{x}, \; \frac{x^p + (1+y^q)^p}{x^p y} \right)
\end {equation}
The second system, coming from extended matrix mutation, is given by a piecewise linear map from $\R^2$ to itself and is closely related to the tropicalization of the first system.
It is given explicitly as
\begin {equation} \label{eq:mu_tropical}
 \trop(x,y) = \left( -x + q[y+p[x]_+]_+, \; -y - p[x]_+ \right),
 \end {equation}
where $[x]_+ := \mathrm{max}(x,0)$.\\

Both systems exhibit behavior which is associated with integrability; namely they each preserve a natural symplectic form, and seem to posess a conserved quantity.
This adds to the growing body of work on cluster algebra theory in dynamical systems, integrability, and mathematical physics which the interested reader can get a glimpse of from various surveys \cite{DiF} \cite{GlickRupel} \cite{HLK} \cite{TYsurvey}.
Our systems share some similarities with other systems from cluster algebra mutation, but we also find new phenomena by allowing integers to be replaced with real numbers.\\

\begin {figure}[h!]
\centering
\includegraphics[scale=0.35]{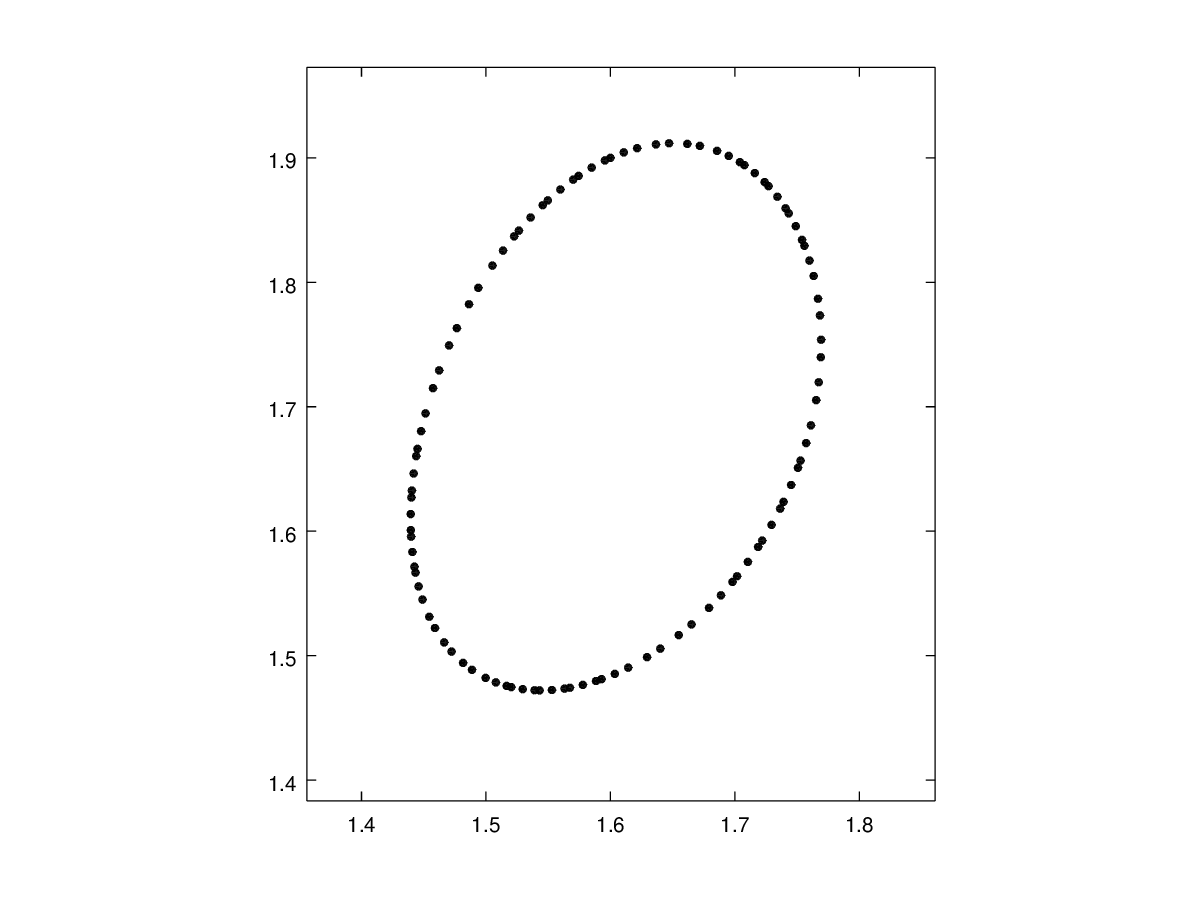}
\caption {An orbit of $\mu$ when $p=1.2734$ and $q=0.8421$ ($100$ iterations pictured)}
\label {fig:orbit}
\end {figure}

We show the orbits of the first system are unbounded when $pq \geq 4$ and conjecture that they are bounded for $pq < 4$.
Furthermore, when $pq < 4$ the orbits empirically appear to lie on simple closed curves and appear aperiodic whenever $p$ and $q$ are not both integers.
An example of a typical orbit is shown in Figure \ref{fig:orbit}.
In the case of cluster mutation (when $p$ and $q$ are integers) the only bounded orbits one obtains come from periodic systems (finite-type cluster algebras).
So, these closed curves for apparently aperiodic systems are a new feature, and it is not clear what a conserved quantity should be.\\

\begin {figure}[h!]
\centering
\includegraphics[scale=0.15]{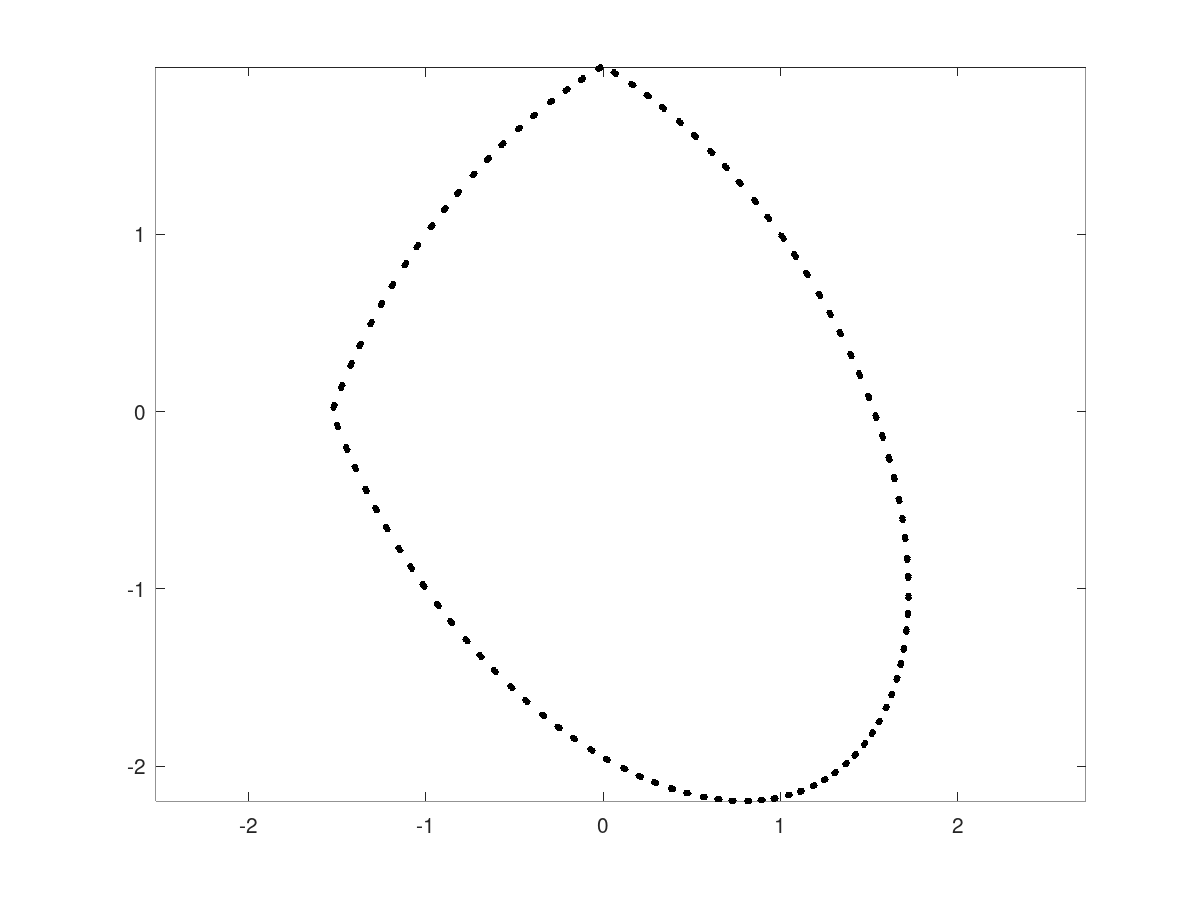}
\caption {An orbit of the tropical map when $p=1.176$ and $q=0.723$ ($200$ iterations pictured)}
\label {fig:tropical_orbit}
\end {figure}

Although it has proven difficult to prove integrability results about the first system, the tropicalized system turns out to be much more accessible.
For the tropicalization, we are able to find a conserved quantity, which establishes a form of integrability. The conserved quantity is a piecewise
function, given by quadratic polynomials on each domain. The level sets are piecewise smooth curves, whose smooth pieces are either ellipses, lines, or hyperbolas,
depending on whether $pq$ is less than, equal to, or greater than $4$, respectively. An example of a typical orbit is pictured in Figure \ref{fig:tropical_orbit}.\\

We are also able to find new examples of periodicity for non-integer values of the parameters $p$ and $q$.
These new instances of periodic orbits occur for systems corresponding to non-crystallographic root systems $I_2(m)$, 
for which the non-tropical system has previously been shown to be nearly periodic~\cite{Lampe}.
In the case when $p$ and $q$ are integers, the tropical system is periodic if and only if the original system is periodic.
So this periodicity for certain tropical systems, without periodicity for the corresponding non-tropical system, is another new trait of the real valued mutation generalization.
\\

In Section \ref{sec:x}, we give a more detailed definition of the first map given above. We then show that it is
symplectic, and that the orbits are unbounded when $pq \geq 4$. \\

In Section \ref{sec:trop}, we give a more detailed definition of the second map given above. The connection between this
map and the tropicalization of the first map is discussed. This map is also shown to be symplectic (for a different symplectic form),
and we give an explicit conserved quantity. Finally, we prove that for certain values of $p$ and $q$, the tropical map is
periodic. We also discuss the connections between our results and the \emph{asymptotic sign-coherence} observed in \cite{ASC}.

\bigskip

\section {The Dynamical System}
\label{sec:x}

Let $p$ and $q$ be two positive real numbers. 
Consider the following two maps on the positive quadrant $\Bbb{R}^+ \times \Bbb{R}^+$:
\[ 
    \mu_1(x,y) = \left( \frac{1+y^q}{x}, \, y \right) \hspace{0.5cm} \text{and} \hspace{0.5cm}
    \mu_2(x,y) = \left( x, \, \frac{1+x^p}{y} \right)
\]

If $p$ and $q$ are integers, then these maps are the mutations for the cluster algebra with the following exchange matrix \cite{fz02}:
\[ B = \begin{pmatrix} 0 & -p \\ q & 0 \end{pmatrix} \]
The exchange matrix for a cluster algebra is required to have integer entries, but here we allow arbitrary real numbers.
Three-by-three exchange matrices with real entries were studied in \cite{FT19, FL} and classification of finite type 
exchange matrices with real entries has been given by Felikson and Tumarkin~\cite{FT}.\\

We denote the composition of both mutations by $\mu := \mu_2 \circ \mu_1$, and it is given by Equation (\ref{eq:mu}).
When $p,q$ are integers, this is the usual case of a cluster algebra from a skew-symmetrizable matrix,
and when $q=1$ and $p \in \{1,2,3\}$ the cluster algebra is of finite type, and hence $\mu$ is periodic \cite{fz03}.
The purpose of the present paper is to investigate the properties of the dynamical system generated by
the map $\mu$ when $p$ and $q$ are arbitrary real numbers.

\bigskip

It is worth mentioning that Lampe studied the special case of these systems when $p=1$ and $q = 4 \cos^2(\pi/m)$
for an integer $m \geq 3$, proving that they are approximately periodic in a quantifiable sense~\cite{Lampe}.

\bigskip

The map $\mu$ exhibits some properties commonly associated with integrability, namely the preservation of a symplectic form. 

\begin {prop}
    The map $\mu$ preserves the following symplectic form: 
    \[ \omega = d \log(x) \wedge d \log(y) = \frac{1}{xy} \, dx \wedge dy \]
    Equivalently, $\mu$ is a Poisson map for the bracket on $\Bbb{R}^+ \times \Bbb{R}^+$ defined by $\{x,y\} = xy$.
\end {prop}
\begin {proof}
    Recall that $\mu_1(x,y) = \left( \frac{1+y^q}{x}, \, y \right)$. 
    Therefore
    \[ \mu_1^* dx = d \left( \frac{1+y^q}{x} \right) = - \frac{1+y^q}{x^2} \, dx + q \, \frac{y^{q-1}}{x} \, dy \]
    Therefore the pullback of $\omega$ by $\mu_1$ is
    \[ \mu_1^* \omega = \frac{x}{y(1+y^q)} \, d \left( \frac{1+y^q}{x} \right) \wedge dy = \frac{-1}{xy} \, dx \wedge dy = -\omega \]
    Similarly, one can see that $\mu_2^* \omega = -\omega$. Therefore the composition gives $\mu^* \omega = \omega$.
\end {proof}

\bigskip

In the two-dimensional setting, a symplectic mapping $\mu$ is integrable (in the Liouville sense) if it possesses a conserved quantity;
i.e. if there is some function $H$ so that $H \circ \mu = H$.
However, we have not been able to find a conserved quantity for this map $\mu$.
On the other hand, we study the tropical version of this map in Section \ref{sec:trop}, and in the tropical setting
we are able to find a conserved quantity.

\bigskip

\subsection {Orbits are Unbounded When $pq \geq 4$} \label{subsec:unbounded}

\bigskip

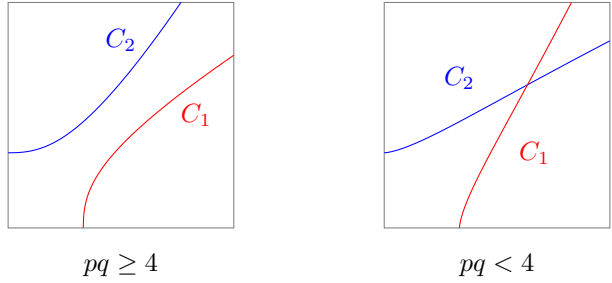
\begin {figure}[h!]
\centering
\begin {tikzpicture}
    \draw[gray] (0,0) -- (3,0) -- (3,3) -- (0,3) -- cycle;
    \draw[smooth, samples=100, domain=0:2.3, variable=\x, blue] plot ({\x}, {sqrt(1+\x^(2.5))});
    \draw[smooth, samples=100, domain=0:2.3, variable=\x, red] plot ({sqrt(1+\x^(2.5))}, {\x});

    \draw [red]  (2.5,1.5) node {$C_1$};
    \draw [blue] (1.5,2.5) node {$C_2$};

    \draw (1.5,-0.5) node {$pq \geq 4$};

    \begin {scope} [shift={(5,0)}]
        \draw[gray] (0,0) -- (3,0) -- (3,3) -- (0,3) -- cycle;
        \draw[smooth, samples=100, domain=0:3, variable=\x, blue] plot ({\x}, {sqrt(1+\x^(1.5))});
        \draw[smooth, samples=100, domain=0:3, variable=\x, red] plot ({sqrt(1+\x^(1.5))}, {\x});

        \draw [red]  (2,1) node {$C_1$};
        \draw [blue] (1,2) node {$C_2$};

        \draw (1.5,-0.5) node {$pq < 4$};
    \end {scope}
\end {tikzpicture}
\caption {Fixed points of the maps $\mu_1$ and $\mu_2$.}
\label {fig:fixed_curves}
\end {figure}

Consider the case where $pq \geq 4$. It is easily seen that the set $C_1$ of fixed points of $\mu_1$ is the curve $x^2 = 1+y^q$ and the
set $C_2$ of fixed points of $\mu_2$ is the curve $y^2=1+x^p$. The typical situation is pictured in Figure \ref{fig:fixed_curves}.

\bigskip

The maps $\mu_1$ and $\mu_2$ interchange the regions on opposite sides of the curves $C_1$ and $C_2$.
This is stated in the following proposition.

\bigskip

\begin {prop}
    Let $(x',y') = \mu_1(x,y)$ and $(x'',y'') = \mu_2(x,y)$.
    \begin {enumerate}
    \item[$(a)$] If $x < \sqrt{1+y^q}$, then $x' > \sqrt{1+y^q}$, and vice versa. 
    \item[$(b)$] If $y < \sqrt{1+x^p}$, then $y'' > \sqrt{1+x^p}$ and vice versa.
    \end {enumerate}
    \label{prop:reflect}
\end {prop}
\begin {proof}
    The formula for $\mu_1$ says that $x' = \frac{1+y^q}{x}$. If we assume that $x < \sqrt{1+y^q}$, then
    \[ x' > \frac{1+y^q}{\sqrt{1+y^q}} = \sqrt{1+y^q} \]
    The other statements are similarly simple.
\end {proof}

\bigskip

\begin {figure}[h]
\centering
\begin {tikzpicture}
    \draw[gray] (0,0) -- (3,0) -- (3,3) -- (0,3) -- cycle;
    \draw[smooth, samples=100, domain=0:2, variable=\x, blue] plot ({\x}, {1+\x});
    \draw[smooth, samples=100, domain=0:1.315, variable=\x, red] plot ({(1+\x^(1.25))^(1.25)}, {\x});

    \draw (1.5,1.5) node {I};
    \draw (1.8,0.3) node {II};
    \draw (0.3,1.8) node {III};

    \draw [red]  (2.5,1.5) node {$C^{u,v}_1$};
    \draw [blue] (1.0,2.5) node {$C^{u,v}_2$};
\end {tikzpicture}
\caption {The regions defined by $C^{u,v}_1$ and $C^{u,v}_2$.}
\label {fig:regions}
\end {figure}
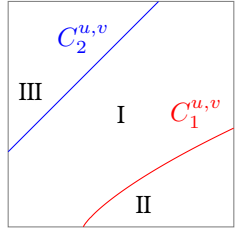

We will look at the map in a different set of coordinates: $(u,v) = (x^p,y^2)$. In these coordinates,
the curve $C_2$ is the line $v=1+u$ and $C_1$ is the curve $u=(1+v^{q/2})^{p/2}$.
We will use the notation $C^{u,v}_1$ and $C^{u,v}_2$ to denote the curves $C_1$ and $C_2$ respectively when we wish to emphasize we are considering the coordinates $(u,v)$. 
The mutation maps in these coordinates look like
\[ \mu_1(u,v) = \left( \frac{(1+v^{q/2})^p}{u}, \; v \right) \quad \text{ and } \quad \mu_2(u,v) = \left( u, \; \frac{(1+u)^2}{v} \right) \]
The curves $C_1^{u,v}$ and $C_2^{u,v}$ divide the first quadrant into three regions. 
In the $(u,v)$-coordinates, the regions are pictured in Figure \ref{fig:regions}.

\bigskip

\begin {lem} \label {lem:increasing}
    Let $H(v)$ be the horizontal distance between $C^{u,v}_1$ and $C^{u,v}_2$ at height $v$, and let $V(u)$ be the vertical distance
    between $C^{u,v}_1$ and $C^{u,v}_2$ as a function of $u$. Then $H$ and $V$ are both increasing functions on the domain $[1,\infty)$.
\end {lem}
\begin {proof}
    The horizontal distance is given by
    \[ H(v) = (1+v^{q/2})^{p/2} - v + 1 \]
    The derivative is
    \[ \frac{dH}{dv} = \frac{pq}{4} \, (1+v^{q/2})^{p/2-1}v^{q/2-1} - 1 \]
    Note that $(1+v^{q/2})^{p/2-1} > v^{(q/2)(p/2-1)}$, and $\frac{pq}{4} \geq 1$, so we have
    \[ \frac{dH}{dv} > v^{(q/2)(p/2-1)} \cdot v^{q/2-1} - 1 = v^{pq/4-1} - 1 \]
    Recall that we only consider $v \geq 1$.
    Since $pq \geq 4$, then $v^{pq/4-1} \geq 1$, and so $\frac{dH}{dv} > 0$.

    \bigskip

    The vertical distance is given by
    \[ V(u) = 1 + u - \left( u^{2/p} - 1 \right)^{2/q} \]
    The derivative is given by
    \[ \frac{dV}{du} = 1 - \frac{4}{pq} \left( u^{2/p} - 1 \right)^{2/q - 1} u^{2/p - 1} \]
    Since $\left( u^{2/p} - 1 \right)^{2/q - 1} > u^{(2/p)(2/q -1)}$, and $\frac{4}{pq} \leq 1$, we have
    \[ \frac{dV}{du} > 1 - u^{(2/p)(2/q-1)} \cdot u^{2/p-1} = 1 - u^{4/pq - 1} \]
    Since $4/pq - 1 < 0$, this means $u^{4/pq - 1} < 1$, and so $\frac{dV}{du} > 0$.
\end {proof}

\bigskip

\begin {lem} \label {lem:unbounded}
    Let $\mu = \mu_2 \circ \mu_1$, and let $(u,v)$ be a point in the first quadrant.
    \begin {itemize}
        \item[$(a)$] If $(u,v)$ is in region III, then the orbit under $\mu$ is unbounded. Moreover, there exists a constant $K > 0$ independent of the point $(u,v)$ such that $u' - u > K$ where $\mu(u,v) = (u', v')$.
        \item[$(b)$] If $(u,v)$ is in region II, then the orbit under $\mu^{-1}$ is unbounded. Moreover, there exists a constant $K > 0$ independent of the point $(u,v)$ such that $v' - v > K$ where $\mu^{-1}(u,v) = (u', v')$.
    \end {itemize}
\end {lem}
\begin {proof}
    Assume $(u,v)$ is in region III. Denote $\mu(u,v) = (u',v')$. Then $\mu_1$ takes $u$ to $u'$ with $u' > u$ since by Proposition~\ref{prop:reflect} it is the case that $(u',v)$ in region II. It follows that $v' > v$ since $(u',v')$ is in region III by Proposition~\ref{prop:reflect} and $C_2^{u,v}$.
    This means that $H(v') > H(v)$ by Lemma \ref{lem:increasing}. The change in the $u$-value under the map $\mu_1$ (i.e. $u'-u$) is at least as large as $H(v)$,
    since $\mu_1$ sends a point in region III to a point in region II (at the same height $v$).

    \medskip

    This means if we define $(u_k,v_k) = \mu^k(u,v)$, then $u_k-u_{k-1} > u_{k-1} - u_{k-2}$. Thus, the
    $u$-coordinate will become arbitrarily large after sufficiently many iterations.
    Now to show the constant $K$ exists we see that if $(u,v)$ is in region III it is that case that $v \geq 1$.
    So, $u' - u \geq H(v) \geq H(1) = K > 0$ gives us our constant.

    \bigskip

    The proof of $(b)$ is the same, up to reversing the roles of $u$ and $v$, and reversing the roles of $\mu_1$ and $\mu_2$.
\end {proof}

\begin {thm}
    If $pq>4$, then all orbits of $\mu = \mu_2 \circ \mu_1$ are unbounded.
    \label{thm:unbounded}
\end {thm}
\begin {proof}
    By Lemma \ref{lem:unbounded}, it suffices to show that for any $(u,v)$, the sequence of iterates $\mu^k(u,v)$ will
    eventually reach region III.
    If the point starts in region III, then there is nothing to show. If the point $(u,v)$ is in region I, then after one
    iteration, $\mu(u,v)$ will be in region III. So it only remains to see that it is impossible for $\mu^k(u,v)$ to be in
    region II for all $k$.

    \bigskip

    Suppose that there is some $(u,v)$ in region II for which $\mu^k(u,v)$ remains in region II for all $k$. Let $R$ denote region II.
    Then this would imply that $(u,v)$ belongs to the intersection $\bigcap_{k = 0}^\infty \mu^{-k}(R)$.
    If we let 
    \[ m_k = \inf_{(u,v) \in \mu^{-k}(R)} v, \] 
    then this implies the sequence $m_k$ increases without bound since the constant $K > 0$ in Lemma~\ref{lem:unbounded} implies that the $v$-coordinate will be at least $1 + (k-1)K$ for every point in  $\mu^{-k}(R)$. Therefore we must have $\bigcap_{k=0}^\infty \mu^{-k}(R) = \varnothing$, which is a contradiction.
\end {proof}

\bigskip

\section {Tropical Dynamics}
\label{sec:trop}

As mentioned earlier, we are not able to find a conserved quantity for the map $\mu$ to establish integrability. In this section,
we study a map which is essentially the tropicalization of $\mu$, and in this tropical setting, we are able to prove both that the
map is symplectic and that it has a conserved quantity. 

\subsection {Matrix Mutation}

Fix $p,q,s,t \in \R$ with $p,q > 0$, and consider the matrices
\[B = \begin{pmatrix} 0&p\\-q&0 \\ s & t\end{pmatrix} \quad \text{and} \quad B' = \begin{pmatrix}0&-p\\q&0\\s&t\end{pmatrix} \]
Matrix mutation was defined in \cite{fz02}. We will not give the full definition here, since we only consider a special case.
In the case of a $3$-by-$2$ matrix, the two mutations are given by
\[ \trop_1(B) = \begin{pmatrix} 0&-p\\q&0\\-s&t+p[s]_+\end{pmatrix} \quad \text{and} \quad \trop_2(B') = \begin{pmatrix}0&p\\-q&0\\s+q[t]_+&-t\end{pmatrix} \]
Recall that $[x]_+$ denotes $\mathrm{max}(x,0)$.\\

The matrix mutations will change the last row $(s,t)$,
giving maps $\Bbb{R}^2 \to \Bbb{R}^2$, which we will also refer to as $\trop_1$ and $\trop_2$. That is, we define
\[ \trop_1(s,t) = (-s, \, t + p[s]_+) \quad \text{ and } \quad \trop_2(s,t) = (s+q[t]_+, \, -t) \]
When $p$ and $q$ are integers, these maps are examples of what Reading~\cite{Reading} refers to as \emph{mutation maps}.

\bigskip

As before, we denote $\trop := \trop_2 \circ \trop_1$ which is the map given in~(\ref{eq:mu_tropical}). 
We refer to this $\trop$ (as opposed to the $\mu$ from the earlier sections) as the \emph{tropical map}, 
and we refer to the dynamical system generated by $\trop$ as the \emph{tropical dynamical system}.

\bigskip

\begin {prop}
    The maps $\trop_1$ and $\trop_2$ satisfy $\trop_1^* \omega = \trop_2^* \omega = -\omega$ for the standard symplectic form $\omega = ds \wedge dt$.
    The composition $\trop = \trop_2 \circ \trop_1$ is therefore symplectic. Equivalently, $\trop$ is a Poisson map for the bracket
    given by $\{s,t\} = 1$.
\end {prop}
\begin {proof}
    The maps $\trop_1$ and $\trop_2$ are piecewise linear, so the form of the maps depends on which region of the domain
    the point $(s,t)$ lies in. For $\trop_1$, we have either $\trop_1(s,t) = (-s,t)$ or $\trop_1(s,t) = (-s,t+ps)$.
    In the first case, clearly $\trop_1^* (ds \wedge dt) = - ds \wedge dt$. In the second case,
    \[ \trop_1^* \omega = d(-s) \wedge d ( t + ps ) = - ds \wedge \left( dt + p ds \right) = - ds \wedge dt \]
    The calculation for $\trop_2$ is similar.
\end {proof}

\bigskip

\subsection {Relation With Tropicalization}

The \emph{tropicalization} of an expression $f(x_1,\dots,x_n)$ is the following limit, provided that it exists:
\[ \hat{f}(x_1,\dots,x_n) = \lim_{t \to \infty} \frac{1}{t} \, \log f \left( e^{tx_1}, \dots, e^{tx_n} \right) \]
Two particular examples of interest are when $f(x,y) = x+y$ and $g(x,y) = xy$. In the former case, the tropicalization
is $\hat{f}(x,y) = \mathrm{max}(x,y)$, and in the latter, $\hat{g}(x,y) = x+y$.

\bigskip

We will discuss in this section how $\trop_1$ and $\trop_2$ are essentially the tropicalizations
of $\mu_1$ and $\mu_2$. The following results can be verified by direct calculation.

\bigskip

\begin {prop}
    The tropicalizations of $\mu_1$ and $\mu_2$ are given by
    \[ \hat{\mu}_1(x,y) = \left( q [y]_+ - x, \, y \right) \quad \text{and} \quad \hat{\mu}_2(x,y) = \left( x, \, p[x]_+ - y \right) \]
\end {prop}

\bigskip

\begin {prop}
    Define the reflections $R_x(x,y) = (-x,y)$ and $R_y(x,y) = (x,-y)$. The extended matrix mutation maps $\trop_k$
    are related to the tropicalized maps $\hat{\mu}_k$ by
    \[ \trop_2 = R_y \circ \hat{\mu}_1 \circ R_x \quad \text{and} \quad \trop_1 = R_x \circ \hat{\mu}_2 \circ R_y \]
    In particular, the composition $\trop = \trop_2 \circ \trop_1$ is given by
    \[ \trop = R_y \circ \hat{\mu}_1 \circ \hat{\mu}_2 \circ R_y \]
    In other words, up to conjugation by $R_y$, the map $\trop$ is a tropical version of $\mu^{-1}$.
\end {prop}

\bigskip

\subsection {A Conserved Quantity}

\bigskip

Let us consider the functions
\[f(x,y) = px^2 + pqxy + qy^2 \quad \text{and} \quad g(x,y) = px^2 - pqxy + qy^2 \]
which will be essential to the study of our tropical dynamical system.
Our map $\trop$ is piecewise linear, and the following lemma shows that $f$ and $g$ are conserved quantities for certain domains of linearity.

\begin{lem}
    For any $(s,t) \in \R^2$ we have that
    \begin {itemize}
        \item[$(a)$] $f(s + qt,-t) = f(s,t) = f(-s, t + ps)$
        \item[$(b)$] $f(s,-t) = g(s,t) = f(-s, t)$
    \end {itemize}
\label{lem:f}
\end{lem}
\begin{proof}
    Part $(b)$ is obvious. We will show the calculation for part $(a)$.
    Let us show that $f(s+qt,-t) = f(s,t)$.
    The fact that $f(-s,t + ps) = f(s,t)$ is completely analogous.
    We perform the calculation:
    \begin{align*}
    f(s+qt,-t) &= p(s+qt)^2 +pq(s+qt)(-t) + q(-t)^2\\
    &= ps^2 + 2pqst + pq^2t^2 -pqst - pq^2t^2 + qt^2\\
    &= ps^2 + pqst + qt^2\\
    &= f(s,t)
    \end{align*}
\end{proof}

\bigskip

We are now ready to give our main theorem on this tropical dynamical system.
Define the following piecewise function using $f$ and $g$:
\[
    \phi(s,t) = \begin{cases} 
        g(s,t) & s < 0 \text{ and } t > 0 \\[1ex]
        f(s,t) & \text{otherwise} 
    \end{cases} 
\]

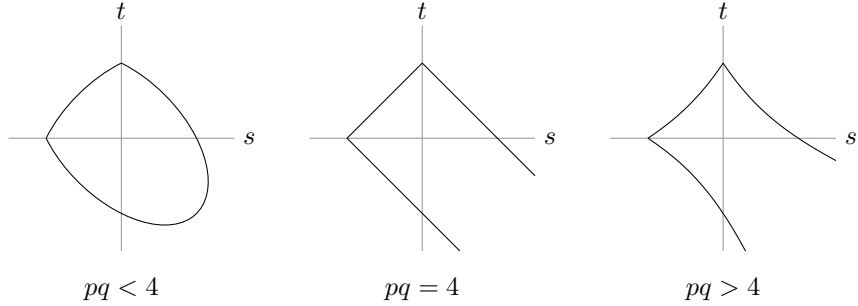
\begin {figure}[h]
\centering
\begin {tikzpicture}
    \draw[gray!80] (-1.5,0) -- (1.5,0);
    \draw[gray!80] (0,-1.5) -- (0,1.5);

    \draw (1.7,0) node {$s$};
    \draw (0,1.7) node {$t$};

    \draw[rotate=-45] ({-1/sqrt(2)},{-1/sqrt(2)}) arc (-120:120:{sqrt(2)} and {sqrt(2/3)});
    \draw[rotate=-45] ({-1/sqrt(2)},{1/sqrt(2)})  arc (150:210:{sqrt(2/3)} and {sqrt(2)});

    \draw (0,-2) node {$pq < 4$};

    \begin {scope} [shift={(4,0)}]
        \draw[gray!80] (-1.5,0) -- (1.5,0);
        \draw[gray!80] (0,-1.5) -- (0,1.5);
        \draw (1.7,0) node {$s$};
        \draw (0,1.7) node {$t$};

        \draw (0,1) -- (1.5,-0.5);
        \draw (-1,0) -- (0.5,-1.5);
        \draw (-1,0) -- (0,1);

        \draw (0,-2) node {$pq = 4$};
    \end {scope}

    \begin {scope} [shift={(8,0)}]
        \draw[gray!80] (-1.5,0) -- (1.5,0);
        \draw[gray!80] (0,-1.5) -- (0,1.5);
        \draw (1.7,0) node {$s$};
        \draw (0,1.7) node {$t$};

        \draw[smooth, domain=0:1.5, variable=\x]  plot ({\x}, {(-3*\x + sqrt(5*\x*\x + 4))/2});
        \draw[smooth, domain=-1:0.3, variable=\x] plot ({\x}, {(-3*\x - sqrt(5*\x*\x + 4))/2});
        \draw[smooth, domain=-1:0, variable=\x] plot ({\x}, {(3*\x + sqrt(5*\x*\x + 4))/2});

        \draw (0,-2) node {$pq > 4$};
    \end {scope}
\end {tikzpicture}
\caption {Level sets for the conserved quantity $\phi$.}
\label {fig:level_sets}
\end {figure}

The level sets of $\phi$ are pictured in Figure \ref{fig:level_sets}. Theorem \ref{thm:conserved} guarantees that the orbits of the tropical
dynamical system stay on these curves.

\bigskip

\begin{thm} \label{thm:conserved}
    The function $\phi(s,t)$ is a conserved quantity for the tropical dynamical system determined by $\trop$.
    In other words, $\phi \circ \trop = \phi$.
\end{thm}
\begin{proof}
    Take any $(s,t) \in \mathbb{R}^2$.
    We let $(s',t') = \trop_1(s,t)$ and $(s'', t'') = \trop_2(s',t') = \trop(s,t)$.\\

    Assume that $s, t \geq 0$, then $(s',t') = (-s, t + ps)$ and $t' \geq 0$.
    So we have $(s'', t'') = (s' + qt', -t')$. Since $t'' \leq 0$, this means $\phi(s'',t'') = f(s'',t'')$.
    By Lemma~\ref{lem:f} we have that $f(s,t) = f(s',t') = f(s'',t'')$ and so 
    \[ \phi(s,t) = f(s,t) = f(s'',t'') = \phi(s'',t'') \]

    \medskip

    Assume that $s \leq 0$ and $t \geq 0$.
    In this case $(s',t') = (-s,t)$ and $g(s,t) = f(s',t')$ by Lemma~\ref{lem:f}.
    Then $(s'',t'') = (s' + qt', -t')$ and $g(s,t) = f(s', t') = f(s'', t'')$ by Lemma~\ref{lem:f}.
    Thus we have $\phi(s,t) = \phi(s'',t'')$ as $t'' \leq 0$.

    \medskip

    Assume that $s \leq 0$ and $t \leq 0$.
    Then $(s',t') = (-s,t)$ and since $t' \leq 0$ it follows $(s'', t'') = (s', -t') = (-s, -t)$.
    Hence, $f(s,t) = f(s'',t'')$ and $\phi(s,t) = \phi(s'',t'')$ as $s'' \geq 0$.\\

    Assume that $s \geq 0$ and $t \leq 0$.
    This means $(s', t') = (-s, t + ps)$ and $f(s,t) = f(s',t')$  by Lemma~\ref{lem:f}.
    If $t+ps \geq 0$ we have that $(s'',t'') = (s' + qt', -t')$ and $f(s,t) = f(s'',t'')$ again by Lemma~\ref{lem:f}.
    In this situation we find $\phi(s,t) = \phi(s'',t'')$ since $t'' \leq 0$.
    Otherwise $t + ps \leq 0$ and $(s'',t'') = (s',-t')$ and so $f(s,t) = g(s'',t'')$ by Lemma~\ref{lem:f}.
    Thus, once again $\phi(s,t) = \phi(s'',t'')$ since $t'' \geq 0$ while $s'' \leq 0$.
\end{proof}

\bigskip

\begin {rmk}
    It is possible to recover the polynomials $f$ and $g$, and to prove Theorem \ref{thm:conserved}, 
    from the results in \cite{bbh_11} and \cite{FT19} concerning the ``\emph{Markov polynomial}''  for rank-3 cluster mutation. 
    More specifically, one may define a skew-symmetric $3 \times 3$ matrix, closely related to $B$, whose Markov constant is
    either $f(s,t) + pq$ or $g(s,t) + pq$, depending on the signs of $s$ and $t$. Since the Markov constant is invariant under
    mutations, and since $pq$ is constant, this implies $f$ and $g$ are conserved (depending again on signs of $s,t$).
\end {rmk}

\bigskip

\subsection {Some Periodic Orbits}

\bigskip

We now wish to prove that for certain choices of $p$ and $q$ our tropical dynamical system is periodic.
To begin we define two linear maps which agree with $\trop_1$ and $\trop_2$ respectively on certain domains.
We let 
\begin{align*}
\tau_1(s,t) &= (-s, t + ps)\\
\tau_2(s,t) &= (s + qt, -t)
\end{align*} 
and $\tau = \tau_2 \circ \tau_1$.
Also, denote the Chebyshev polynomials of the second kind by $U_n(x)$ which can be defined by the recurrence
\[U_n(x) = 2xU_{n-1}(x) - U_{n-2}(x)\]
with $U_0(x) = 1$ and $U_1(x) = 2x$.
These polynomials satisfy 
\begin{equation}
U_n(\cos(\theta)) = \frac{\sin((n+1)\theta)}{\sin(\theta)}
\label{eq:cheb}
\end{equation}
and interact nicely with $\tau$.
Let us establish the notation
\begin{align*}
\kappa &= \sqrt{pq}  &\nu &= \sqrt{\frac{p}{q}}
\end{align*}
to denote the square roots of the product and ratio of $p$ and $q$.
We now have the following lemma which describes the action of $\tau$.

\begin{lem}[{\cite[Proposition 3.1]{ASC}}]
For any $(s,t) \in \R^2$ let $(s_n, t_n) = \tau^n (s,t)$ and $(\tilde{s}_n, \tilde{t}_n) = \tau_1 \circ \tau^n (s,t)$.
Then for any $n$ we have
\begin{align*}
s_n &= s U_{2n}\left(\frac{\kappa}{2}\right) + t \nu^{-1} U_{2n-1}\left(\frac{\kappa}{2}\right)\\
t_n &= -s \nu  U_{2n-1} \left(\frac{\kappa}{2}\right) - t U_{2n-2}\left(\frac{\kappa}{2}\right)\\
\tilde{s}_n &= -s U_{2n} \left(\frac{\kappa}{2}\right) - t \nu^{-1} U_{2n-1}\left(\frac{\kappa}{2}\right)\\
\tilde{t}_n &= s \nu U_{2n+1} \left(\frac{\kappa}{2}\right) + t U_{2n}\left(\frac{\kappa}{2}\right)
\end{align*}
\label{lem:tau}
\end{lem}
\begin{proof}
The proof follows what is done in~\cite[Proposition 3.1]{ASC} based on recurrence of Chebyshev polynomials of the second kind.
The only difference here is different convention on initial point; so, the indices are shifted.
Also, since we are not assuming $pq \geq 4$ we use $\tau$ instead of $\mu$ as to not have to worry if $s_n$ or $t_n$ is positive.
\end{proof}

When $0< pq < 4$ we will use the notation $pq = 4 \cos^2(\theta)$ for $0 < \theta < \frac{\pi}{2}$.
Since the cosine function maps the open interval $(0, \frac{\pi}{2})$ bijectively to the open interval $(0,1)$ its square does as well and there is a unique such $\theta$.
It follows that $\frac{\kappa}{2} = \cos(\theta)$ or equivalently that $\kappa = 2 \cos(\theta)$.

\begin{lem}
Assume $pq < 4$.
For any $(s,t) \in \R^2$ let $(s_n, t_n) = \tau^n (s,t)$ and $(\tilde{s}_n, \tilde{t}_n) = \tau_1 \circ \tau^n (s,t)$.
Then for any $n$ we have
\begin{align*}
s_n &=\frac{s\sin((2n+1)\theta)}{\sin(\theta)} +\frac{t\sin(2n\theta)}{\nu\sin(\theta)}\\
t_n &= -\frac{s\nu\sin(2n\theta)}{\sin(\theta) } - \frac{t\sin((2n-1)\theta)}{\sin(\theta)}\\
\tilde{s}_n &= -\frac{s\sin((2n+1)\theta)}{\sin(\theta)} -\frac{t\sin(2n\theta)}{\nu\sin(\theta)}\\
\tilde{t}_n &= \frac{s\nu\sin((2n+2)\theta)}{\sin(\theta)} + \frac{t\sin((2n+1)\theta)}{\sin(\theta)}
\end{align*}
\label{lem:tautheta}
\end{lem}
\begin{proof}
This follows from Lemma~\ref{lem:tau} and Equation~(\ref{eq:cheb}) by substitution of $\cos(\theta)$ into the Chebyshev polynomials of the second kind.
\end{proof}

We now find the use of the functions $\tau_1$, $\tau_2$, and $\tau$.
The iteration of $\tau$ turns out to agree with the iteration of $\trop$ for a certain number of iterations.
This allows us to use properties of the Chebyshev polynomials of the second kind in analyzing our dynamical system coming from cluster mutation.

\begin{lem}
Assume $pq < 4$ and $(2n+2)\theta \leq \pi$.
For any $(s,t) \in \R^2$ with $s,t \geq 0$ we have $\trop^{n+1}(s,t) = \tau^{n+1}(s,t)$.
\label{lem:leq}
\end{lem}
\begin{proof}
We will have $\trop^{n+1}(s,t) = \tau^{n+1}(s,t)$ as long as $s_n > 0$ and $\tilde{t}_n > 0$.
Now use Lemma~\ref{lem:tautheta}.
When $(2n+1)\theta \leq \pi$ we will have $s_n > 0$.
Also if $(2n+2)\theta \leq \pi$, then $\tilde{t}_n > 0$.
\end{proof}

It will be particularly important to understand the action of the system when it maps the fourth quadrant to itself.
The following lemmas aid in this understanding. 

\begin{lem}
For $(s,t) \in \R^2$ with $s >0$ and $t < 0$ let $(s',t') = \tau(s,t) = \left((pq-1)s + qt, -ps-t\right)$ and assume $s' > 0$, then
\[\frac{t}{s} - \frac{t'}{s'} = \frac{ps^2 + pqst + qt^2}{ss'} > 0\]
and
\[\arctan\left( \frac{ps^2 + pqst + qt^2}{ss'+tt'} \right)\]
is the angle between the two lines through the origin containing $(s,t)$ and $(s',t')$ respectively.
\label{lem:slope}
\end{lem}
\begin{proof}
We start by computing the change in slope
\begin{align*}
\frac{t}{s} - \frac{t'}{s'}  &= \frac{ts' - st'}{ss'}\\
	&= \frac{t((pq-1)s + qt) - s(-ps - t))}{ss'}\\
	&= \frac{ps^2 + pqst + qt^2}{ss'}
\end{align*}
and we have the formula in the lemma.
Now we have assumed $s'>0$ and the expression $ps^2 + pqst + qt^2$ is positive for any real numbers $s$ and $t$ not both zero.
It follows that the ratio is positive.
Now computing the change in polar angle we find
\begin{align*}
\arctan\left( \frac{t}{s} \right) - \arctan\left( \frac{t'}{s'} \right) &= \arctan\left( \frac{\frac{t}{s} - \frac{t'}{s'}}{1 + \frac{tt'}{ss'}}\right)\\
&=  \arctan\left( \frac{\frac{ps^2 + pqst + qt^2}{ss'}}{\frac{ss'+ tt'}{ss'}}\right)\\
&= \arctan\left( \frac{ps^2 + pqst + qt^2}{ss'+tt'} \right)
\end{align*}
and the lemma is proven.
\end{proof}

Before proving the theorem we give a lemma which reduces the problem of showing periodicity to checking periodicity for initial points in the first quadrant.

\begin{lem}
Fix $c \in \R$ with $c > 0$.
If for all $p,q$ such that $pq = c$ we have that $\trop^k(s,t) = (s,t)$ whenever $(s,t)$ is in the first quadrant, then $\trop^k(s,t) = (s,t)$ for all $(s,t) \in \R^2$ .
\label{lem:Q1}
\end{lem}
\begin{proof}
First observe that $(\trop_1 \circ \trop_2)^k(s,t) = (s,t)$ whenever $(s,t)$ is in the first quadrant since we may exchange both $s$ with $t$ and $p$ with $q$ which does not affect the product $pq$.
Next we note that for any $(s,t)$ in the third quadrant $\trop(s,t) = (-s,-t)$ is in the first quadrant.
Hence, $\trop^{k+1}(s,t) = \trop(s,t)$ which implies $\trop^k(s,t) = (s,t)$ for $(s,t)$ in the third quadrant since $\trop$ is invertible.
For $(s,t)$ in the second quadrant we have that $\trop_1(s,t) = (-s,t)$ is in the first quadrant. 
Then $(\trop_1 \circ \trop_2)^k(-s,t) = (-s,t)$ so
\[(\trop_1 \circ \trop_2)^k(\trop_1(s,t)) = \trop_1(\trop^k(s,t)) = \trop_1(s,t)\]
which implies $\trop^k(s,t) = (s,t)$ for $(s,t)$ in the second quadrant since $\trop_1$ is invertible.\\

It remains to handle the case of $(s,t)$ in the fourth quadrant.
To do this we show that for any $(s,t)$ in the fourth quadrant $\trop^j(s,t)$ is not in the fourth quadrant for some $j$.
This implies $\trop^{k+j}(s,t) = \trop^j(s,t)$ since we have proven the lemma for all other quadrants.
This in turn means $\trop^k(s,t) = (s,t)$ by the invertibility of $\trop$.
Assume $(s,t)$ and $(s',t') = \trop(s,t)$ are both in the fourth quadrant, then by Lemma~\ref{lem:slope} we know that the difference in angle when moving from $(s,t)$ and $(s',t')$ is $\arctan\left( \frac{ps^2 + pqst + qt^2}{ss'+tt'} \right)$.
We see that $ps^2 + pqst + qt^2$ is the conserved quantity from Theorem~\ref{thm:conserved} and $ss' + tt'$ is bounded from above because the graph of the conserved quantity in quadrant four is an ellipse.
So, the change in angle with each iteration of $\trop$ is bounded below by a positive constant.
Hence, the orbit must leave the fourth quadrant and the proof is complete.
\end{proof}

\begin{figure}
\centering
\begin{minipage}[b]{0.45\linewidth}
\includegraphics[scale=0.40]{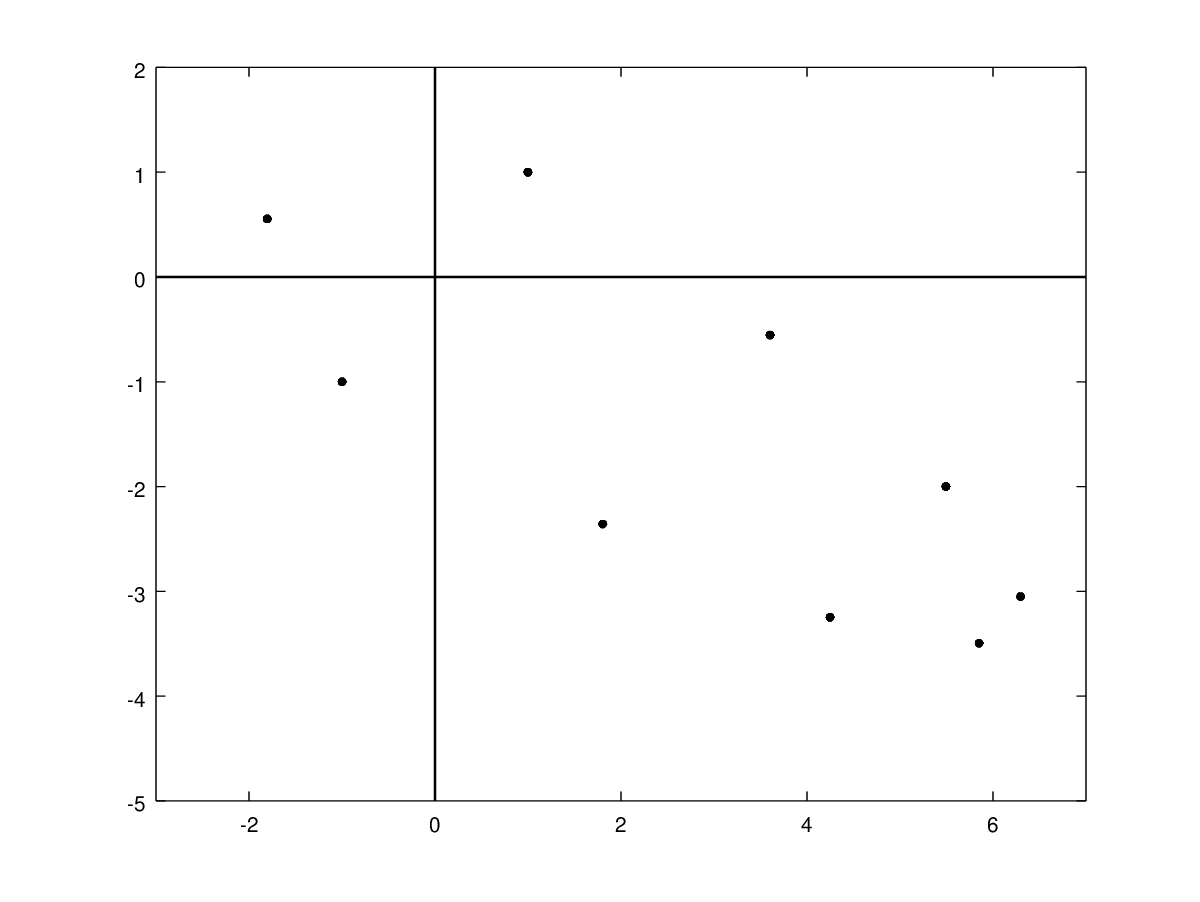}
\end{minipage}
\quad
\begin{minipage}[b]{0.45\linewidth}
\includegraphics[scale=0.40]{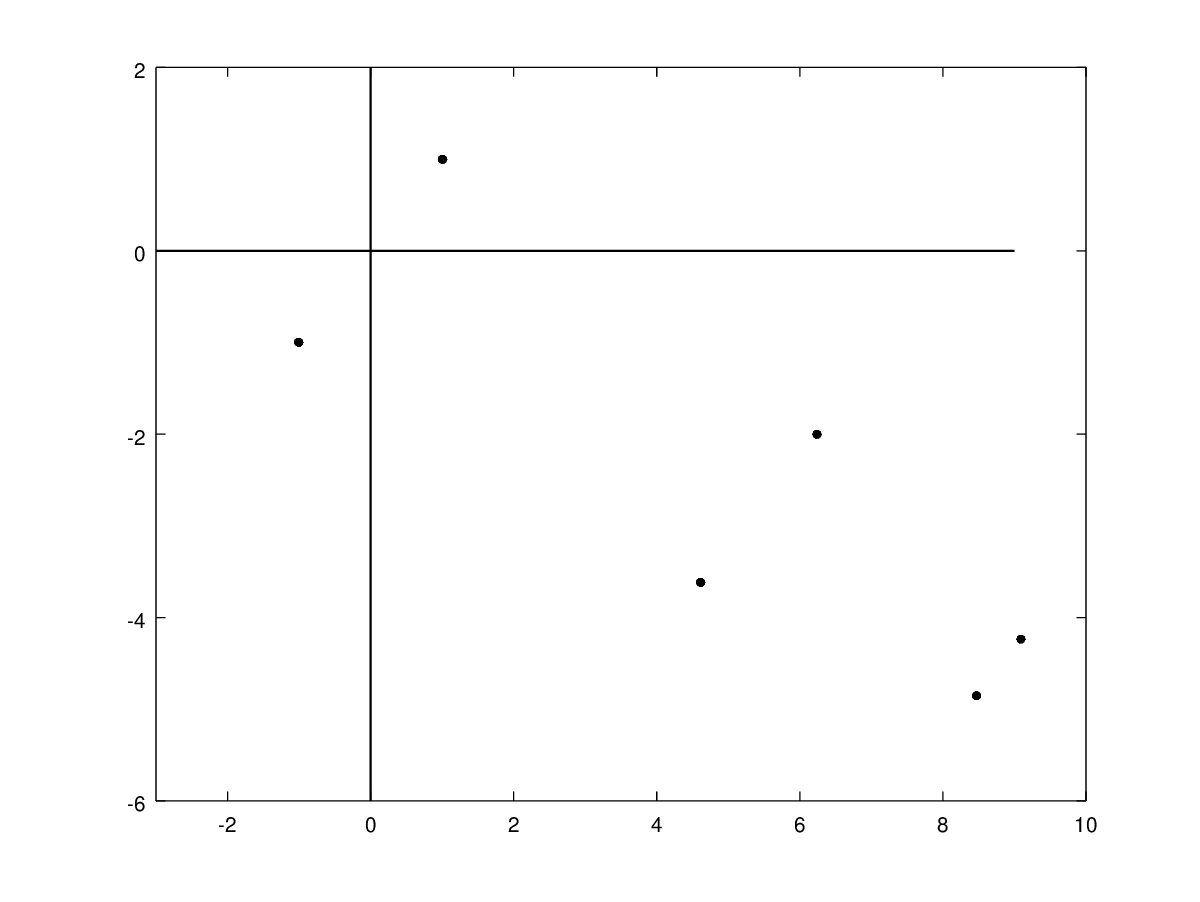}
\end{minipage}
\caption{Example orbits for Theorem~\ref{thm:periodic} with $m=7$ on the left and $m=10$ on the right}
\label{fig:periodic}
\end{figure}

We now prove a theorem on periodic orbits which makes use of the map $\tau$ along with some further analysis.
Two examples of the periodic orbits are shown in Figure~\ref{fig:periodic}.
Let us point out the next theorem applies to what can be considered exchange matrices for non-crystallographic root systems including the matrices
\[
\begin{pmatrix}
0 & 1 \\[1ex] -4\cos^2\left(\frac{\pi}{m}\right) &0
\end{pmatrix}
\]
which were shown by Lampe~\cite{Lampe} to give rise to an almost periodic dynamical system from cluster $x$-variable mutation that we considered in Section~\ref{sec:x}.

\begin{thm}
If $pq = 4\cos^2\left(\frac{\pi}{m}\right)$ with $m \in \mathbb{Z}$ with $m \geq 3$, then $\trop$ is periodic and its orbit has size $m+2$ when $m$ is odd and $\frac{m+2}{2}$ when $m$ is even.
\label{thm:periodic}
\end{thm}

\begin{proof}
Let us denote $\theta = \frac{\pi}{m}$ and assume $(s,t)$ is in the first quadrant.
By Lemma~\ref{lem:leq} we have that $\trop^{n+1}(s,t) = \tau^{n+1}(s,t)$ for $(2n+2) \leq m$.
First assume that $m = 2k$ so that $(2n+2) \leq m$ means $n \leq k - 1$.
So by Lemma~\ref{lem:tautheta} we have
\begin{align*}
\trop^{k}(s, t) &= \tau^{k}(s,t)\\
&= \left(\frac{s\sin((2k+1)\theta)}{\sin(\theta)} +\frac{t\sin(2k\theta)}{\nu\sin(\theta)},  -\frac{s\nu\sin(2k\theta)}{\sin(\theta)}  - \frac{t\sin((2k-1)\theta)}{\sin(\theta)}\right)\\
&= \left(\frac{s\sin(\pi+\theta)}{\sin(\theta)} +\frac{t\sin(\pi)}{\nu\sin(\theta)},  -\frac{s\nu\sin(\pi)}{\sin(\theta)}  - \frac{t\sin(\pi-\theta)}{\sin(\theta)}\right)\\
&= (-s,-t).
\end{align*}
Then since
\[\trop(-s,-t) = \trop_2(\trop_1(-s,-t)) = \trop_2(s,-t) = (s,t)\]
we have that $\trop^{k+1}(s,t) = (s,t)$ and $k+1 = \frac{m}{2} + 1$.\\

Next we assume that $m = 2k+1$ is odd.
Here $(2n+2) \leq m$ means that $n \leq k - \frac{1}{2}$.
Hence, again $n \leq k-1$ since $n$ and $k$ are integers.
By Lemma~\ref{lem:tautheta} we have
\begin{align*}
\trop^{k}(s, t) &= \left(\frac{s\sin((2k+1)\theta)}{\sin(\theta)} +\frac{t\sin(2k\theta)}{\nu\sin(\theta)},  -\frac{s\nu\sin(2k\theta)}{\sin(\theta)}  - \frac{t\sin((2k-1)\theta)}{\sin(\theta)}\right)\\
&= \left(\frac{s\sin(\pi)}{\sin(\theta)} +\frac{t\sin(\pi-\theta)}{\nu\sin(\theta)},  -\frac{s\nu\sin(\pi-\theta)}{\sin(\theta)}  - \frac{t\sin(\pi-2\theta)}{\sin(\theta)}\right)\\
&= \left(\frac{t}{\nu}, -s\nu - 2t\cos(\theta) \right).
\end{align*}
Next we compute
\[\trop_1\left(\frac{t}{\nu}, -s\nu - 2t\cos(\theta) \right) = \left(\frac{-t}{\nu}, -s\nu - 2\cos(\theta)t + \kappa t \right) = \left(\frac{-t}{\nu}, -s\nu \right)\]
and
\[\trop_2 \left(\frac{-t}{\nu}, -s\nu \right) =  \left(\frac{-t}{\nu}, s\nu \right).\]
Hence,
\[\trop^{k+1}(s,t) =  \left(\frac{-t}{\nu}, s\nu \right).\]
Now we find that
\begin{align*}
\trop_1\left(\frac{-t}{\nu}, s\nu \right) &= \left(\frac{t}{\nu}, s\nu \right) & \trop_2\left(\frac{t}{\nu}, -s\nu \right) &= \left(\frac{t}{\nu} + qs\nu, s\nu \right) = \left(\frac{t}{\nu} + s\kappa, -s\nu \right)
\end{align*}
and so
\[\trop^{k+2}(s,t) = \left(\frac{t}{\nu} + s\kappa, -s\nu \right).\]
We let $\hat{s} = \frac{t}{\nu} + s\kappa$ and $\hat{t} = -s\nu$.
Now we will use Lemma~\ref{lem:tautheta} starting from the point $(\hat{s}, \hat{t})$.
Letting $(s_n, t_n) = \tau^n (\hat{s},\hat{t})$ and $(\tilde{s}_n, \tilde{t}_n) = \tau_1 \circ \tau^n (\hat{s},\hat{t})$ we have that
\begin{align*}
s_n &=\frac{\hat{s}\sin((2n+1)\theta)}{\sin(\theta)} +\frac{\hat{t}\sin(2n\theta)}{\nu\sin(\theta)} = \frac{t\sin((2n+1)\theta)}{\nu\sin(\theta)} + \frac{s\kappa\sin((2n+1)\theta)}{\sin(\theta)} - \frac{s\sin(2n\theta)}{\sin(\theta)} \\
t_n &= -\frac{\hat{s}\nu\sin(2n\theta)}{\sin(\theta) } - \frac{\hat{t}\sin((2n-1)\theta)}{\sin(\theta)}\\
\tilde{s}_n &= -\frac{\hat{s}\sin((2n+1)\theta)}{\sin(\theta)} -\frac{\hat{t}\sin(2n\theta)}{\nu\sin(\theta)}\\
\tilde{t}_n &= \frac{\hat{s}\nu\sin((2n+2)\theta)}{\sin(\theta)} + \frac{\hat{t}\sin((2n+1)\theta)}{\sin(\theta)} = \frac{t\sin((2n+2)\theta)}{\sin(\theta)} + \frac{s\kappa\sin((2n+2)\theta)}{\sin(\theta)} -  \frac{s \nu \sin((2n+1)\theta)}{\sin(\theta)}
\end{align*}
and we will further analyze $s_n$ and $\tilde{t}_n$ to show that these values also agree with the map $\trop$ up to a certain value of $n$.
We see that
\begin{align*}
s_n &= \frac{t\sin((2n+1)\theta)}{\nu\sin(\theta)} + \frac{s\kappa\sin((2n+1)\theta)}{\sin(\theta)} - \frac{s\sin(2n\theta)}{\sin(\theta)} \\
&= s \left(\frac{\kappa\sin((2n+1)\theta) - \sin(2n\theta)}{\sin(\theta)}\right) + t \left(\frac{\sin((2n+1)\theta)}{\nu\sin(\theta)}\right)\\
&= s \left(\frac{2\cos(\theta)\sin((2n+1)\theta) - \sin(2n\theta)}{\sin(\theta)}\right) + t \left(\frac{\sin((2n+1)\theta)}{\nu\sin(\theta)}\right)\\
&= s \left(\frac{(\sin((2n+2)\theta) + \sin(2n\theta)) - \sin(2n\theta)}{\sin(\theta)}\right) + t \left(\frac{\sin((2n+1)\theta)}{\nu\sin(\theta)}\right)\\
&= s \left(\frac{\sin((2n+2)\theta)}{\sin(\theta)}\right) + t \left(\frac{\sin((2n+1)\theta)}{\nu\sin(\theta)}\right)
\end{align*}
which is postive for $2n+2 \leq m$.
Next we look at
\begin{align*}
\tilde{t}_n &= \frac{t\sin((2n+2)\theta)}{\sin(\theta)} + \frac{s\kappa\nu\sin((2n+2)\theta)}{\sin(\theta)} -  \frac{s \nu \sin((2n+1)\theta)}{\sin(\theta)}\\
&= s\left(\frac{\kappa\nu\sin((2n+2)\theta) - \nu \sin((2n+1)\theta)}{\sin(\theta)}\right) +  t\left(\frac{\sin((2n+2)\theta)}{\sin(\theta)}\right) \\
&= s\left(\frac{2\cos(\theta)\nu\sin((2n+2)\theta) - \nu \sin((2n+1)\theta)}{\sin(\theta)}\right) +  t\left(\frac{\sin((2n+2)\theta)}{\sin(\theta)}\right) \\
&= s\left(\frac{\nu(\sin((2n+3)\theta) + \sin((2n+1)\theta))  - \nu \sin((2n+1)\theta)}{\sin(\theta)}\right) +  t\left(\frac{\sin((2n+2)\theta)}{\sin(\theta)}\right) \\
&= s\left(\frac{\nu\sin((2n+3)\theta)}{\sin(\theta)}\right) +  t\left(\frac{\sin((2n+2)\theta)}{\sin(\theta)}\right) \\
\end{align*}
which is postive for $2n+3 \leq m$.
Hence, $\tau^{n+1}(\hat{s}, \hat{t}) = \trop^{n+1}(\hat{s}, \hat{t})$ for $2n+3 \leq 2k+1$ or equivalently $n \leq k -1$.
This means
\begin{align*}
\trop^{(k+2) + k}(s,t) &= \trop^{k}(\hat{s}, \hat{t})\\
&= \left(\frac{\hat{s}\sin((2k+1)\theta)}{\sin(\theta)} +\frac{\hat{t}\sin(2k\theta)}{\nu\sin(\theta)},  -\frac{\hat{s}\nu\sin(2k\theta)}{\sin(\theta) } - \frac{\hat{t}\sin((2k-1)\theta)}{\sin(\theta)}\right)\\
&= \left(\frac{\hat{t}}{\nu}, -\nu \hat{s} - \frac{\hat{t} \sin(2\theta)}{\sin(\theta)} \right)\\
&= \left(-s, -t - \nu\kappa s - \frac{2\hat{t} \sin(\theta)\cos(\theta)}{\sin(\theta)} \right)\\
&= (-s, -t - \nu\kappa s - \kappa \hat{t})\\
&= (-s,-t)
\end{align*}
from which it follows that $\trop^{m+2}(s,t) = \trop^{2k+3}(s,t) = \trop(-s,-t) = (s,t)$ and the proof is complete for $(s,t)$ in the first quadrant.
Lemma~\ref{lem:Q1} finishes the proof of the theorem.
\end{proof}

\begin{rmk}
Theorem~\ref{thm:periodic} implies that for
\[ B = \begin{pmatrix} 0 & -p \\ q & 0 \\ a_1 & b_1 \\ \vdots & \vdots \\ a_{\ell} & b_{\ell} \end{pmatrix} \]
with $pq = 4\cos^2\left(\frac{\pi}{m}\right)$ the mutation class of $B$ is finite for any choice of $(a_i, b_i) \in \R^n$ for $1 \leq i \leq \ell$.
In the case of cluster algebras where $p,q,a_i, b_i \in \Z$ all these matrix mutation classes being finite imply that the cluster algebra is of finite type~\cite[Proposition 4.9]{FZ4}.
This would mean the corresponding non-tropical dynamical system would be periodic.
Furthermore, when dealing with cluster algebras (of any rank) periodicity of the tropical system forces periodicity of the non-tropical system with the same period~\cite{Synch}.
However, for $pq = 4\cos^2\left(\frac{\pi}{m}\right)$ when $p$ and $q$ are not both integers, the non-tropical dynamical system does not appear to be periodic, and we can easily check for small $m$ that the non-tropical system does not have the period from Theorem~\ref{thm:periodic}.
\end{rmk}

\subsection{Monotonicity and asymptotic sign coherence}

\bigskip

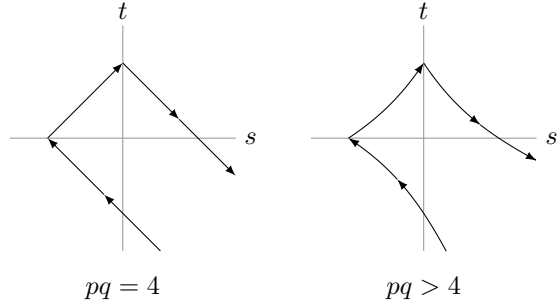
\begin {figure}[h]
\centering
\begin {tikzpicture}
    \begin {scope} [shift={(4,0)}]
        \draw[gray!80] (-1.5,0) -- (1.5,0);
        \draw[gray!80] (0,-1.5) -- (0,1.5);

        \draw (1.7,0) node {$s$};
        \draw (0,1.7) node {$t$};

	\draw[-{latex}] (0.75, 0.25) -- (1.5,-0.5);
        \draw[-{latex}] (0,1) --  (0.75,0.25);
         \draw[-{latex}] (-0.25,-0.75) -- (-1,0);

        \draw[-{latex}] (0.5,-1.5) -- (-0.25,-0.75);

        \draw[-{latex}] (-1,0) -- (0,1);

        \draw (0,-2) node {$pq = 4$};
    \end {scope}

    \begin {scope} [shift={(8,0)}]
        \draw[gray!80] (-1.5,0) -- (1.5,0);
        \draw[gray!80] (0,-1.5) -- (0,1.5);

        \draw (1.7,0) node {$s$};
        \draw (0,1.7) node {$t$};

        \draw[smooth, domain=0:0.75, variable=\x,-{latex}]  plot ({\x}, {(-3*\x + sqrt(5*\x*\x + 4))/2});
         \draw[smooth, domain=0.75:1.5, variable=\x,-{latex}]  plot ({\x}, {(-3*\x + sqrt(5*\x*\x + 4))/2});

        \draw[smooth, domain=-1:-0.35, variable=\x,{latex}-] plot ({\x}, {(-3*\x - sqrt(5*\x*\x + 4))/2});
        \draw[smooth, domain=-0.35:0.3, variable=\x,{latex}-] plot ({\x}, {(-3*\x - sqrt(5*\x*\x + 4))/2});

        \draw[smooth, domain=-1:0, variable=\x,-{latex}] plot ({\x}, {(3*\x + sqrt(5*\x*\x + 4))/2});

        \draw (0,-2) node {$pq > 4$};
    \end {scope}
\end {tikzpicture}
\caption {Level sets for the conserved quantity with direction of motion indicated by arrows.}
\label {fig:motion}
\end {figure}

We will now focus on $pq \geq 4$ and show that the iteration of $\trop$ is monotonic in the sense that it moves in a single direction along the level sets of the conserved quantity.
The direction of motion is shown in Figure~\ref{fig:motion}.
We also show that the orbit of $\trop$ is unbounded and deduce another proof of Gekhtman and Nakanishi's asymptotic sign coherence for rank two cluster algebras~\cite{ASC}.
To begin we collect some inequalities which will be useful to us.

\begin{lem}
Consider $(s,t) \in \R^2$ and let $(s',t') = ((pq-1)s + qt, -ps - t)$.
If $\sqrt{p}s + \sqrt{q} t \geq 0$ and $pq \geq 4$, then $\sqrt{p}s' + \sqrt{q} t' \geq 0$
\label{lem:line}
\end{lem}
\begin{proof}
We compute
\begin{align*}
\sqrt{p}s' + \sqrt{q} t' &= \sqrt{p}pqs - \sqrt{p}s + \sqrt{p}q t - \sqrt{q}ps - \sqrt{q}t\\
&=\sqrt{pq} p \sqrt{q}s - p \sqrt{q} s + \sqrt{p}q t - \sqrt{p}s - \sqrt{q}t\\
&= (\sqrt{pq} - 1) p \sqrt{q} s + \sqrt{p} q t - (\sqrt{p}s + \sqrt{q} t)\\
&\geq  p \sqrt{q} s + \sqrt{p} q t - (\sqrt{p}s + \sqrt{q} t)\\
&= \sqrt{pq} (\sqrt{p}s + \sqrt{q} t) - (\sqrt{p}s + \sqrt{q} t)\\
&= (\sqrt{pq} - 1)(\sqrt{p}s + \sqrt{q} t)\\
&\geq \sqrt{p}s + \sqrt{q} t\\
&\geq 0
\end{align*}
and the lemma is proven.
\end{proof}

\begin{lem}
If $pq \geq 1$, $\sqrt{p}s + \sqrt{q}t \geq 0$ and $t \leq 0$, then $ps + t \geq 0$. 
\label{lem:firstmu}
\end{lem}
\begin{proof}
Multiplying by $\sqrt{p}$ we find that $ps + \sqrt{pq} t \geq 0$.
Since $t \leq 0$ it follows $t \geq \sqrt{pq}t$ as $\sqrt{pq} \geq 1$.
Hence,
\[ps + t \geq ps + \sqrt{pq} t \geq 0\]
and the lemma is proven.
\end{proof}

\begin{lem}
If $pq \geq 1$, $s \leq 0$, and $t \geq 0$, then $\sqrt{p} s' + \sqrt{q} t' \geq 0$ where $s' = -s + qt$ and $t' = - t$.
\label{lem:Q2}
\end{lem}
\begin{proof}
We see that
\[\sqrt{p}(-s + qt) -\sqrt{q}t = - \sqrt{p}s + (\sqrt{pq} -1) \sqrt{q} t \geq 0\]
which holds because $s \leq 0$, $t \geq 0$ and $\sqrt{pq} \geq 1$.
\end{proof}

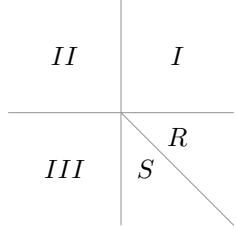
\begin {figure}
\centering
\begin {tikzpicture}
        \draw[gray!80] (-1.5,0) -- (1.5,0);
        \draw[gray!80] (0,-1.5) -- (0,1.5);
        \draw[gray!80] (0,0) -- (1.5, -1.5);

        \draw (0.75,-0.33) node {$R$};
        \draw (0.33,-0.75) node {$S$};
        \draw (0.75,0.75) node {$I$};
         \draw (-0.75,0.75) node {$II$};
        \draw (-0.75,-0.75) node {$III$};
\end{tikzpicture}
\caption{Depiction of regions used in proof of Theorem~\ref{thm:geq4} drawn for case of $p=q$.}
\label{fig:RS}
\end{figure}

We let $\sign: \R \to \{-,0,+\}$ be given by
\[\sign(t) = \begin{cases} + & t > 0 \\ 0 & t = 0 \\ - & t < 0 \end{cases}\]
and for $(s,t) \in \R^2$ let $\sign((s,t)) = (\sign(s), \sign(t))$.

\begin{thm}
Let $(s_0,t_0) \in \R^2 \setminus \{(0,0)\}$, $(s_n, t_n) = \trop^n(s_0,t_0)$, and 
\[\theta_n \in \left( \arctan \left( \frac{-\sqrt{q}}{\sqrt{p}} \right), 2\pi + \arctan \left( \frac{-\sqrt{q}}{\sqrt{p}} \right)\right]\]
 be the angle of $(s_n,t_n)$ in polar coordinates.
If $pq \geq 4$, then
\begin{enumerate}
\item[(i)] the sequence $\{\theta_n\}_{n \geq 0}$ is monotonically decreasing,
\item[(ii)] there exists $N$ such that $\sign((s_n, t_n)) = (+,-)$ for all $n \geq N$,
\item[(iii)] and the sequence $\{(s_n,t_n)\}_{n \geq 0}$ is unbounded.
\end{enumerate}
\label{thm:geq4}
\end{thm}
\begin{proof}
Let $R$ denote the region in the fourth quadrant above the line with slope $\arctan \left( \frac{-\sqrt{p}}{\sqrt{q}} \right)$ and let $S$ denote the region in the fourth quadrant below the line with slope $\arctan \left( \frac{-\sqrt{p}}{\sqrt{q}} \right)$.
In Figure~\ref{fig:RS} we show the regions $R$ and $S$ along with the other quadrants labeled.
To prove (i) and see that $\{\theta_n\}_{n \geq 0}$ is monotonically decreasing we first note that if $(s,t)$ is in the third quadrant then $\trop(s,t) = (-s,-t)$ is in the first quadrant and the polar angle has decreased.
If $(s,t)$ is in the second quadrant then $\trop(s,t) = (-s + qt, -t)$, and this case is covered by Lemma~\ref{lem:Q2} which says that $\trop(s,t)$ will be in region $R$.\\

The case when $(s,t)$ is in the first quadrant and region $R$ are both handled by Lemma~\ref{lem:line}.
For these cases $\trop_1(s,t) = (-s, ps+t)$ and lemma~\ref{lem:firstmu} ensures that $ps+t \geq 0$ so that $\trop(s,t) = ((pq-1)s + qt, -ps - t)$.
So, in this case $\trop(s,t)$ is in the region $R$.
Thus it is clear the polar angle has decreased if $(s,t)$ was in the first quadrant, and if $(s,t)$ was in the fourth quadrant that angle must also decrease by Lemma~\ref{lem:slope}.\\

The last case to prove (i) is when $(s,t)$ is in the region $S$.
In this situation $\trop(s,t)$ can be in quadrant two, three, or four.
If $\trop(s,t)$ is in quadrant two or three the polar angle has certainly decreased.
If $\trop(s,t)$ is in quadrant four then by  Lemma~\ref{lem:slope} the polar angle has decreased.
Thus (i) is proven as we have shown any application of $\trop$ decreases the polar angled within our chosen interval.\\

In order to prove (ii) and (iii) we will start by arguing that by iterating $\trop$ we will always arrive in the region $R$ at some step.
We have already seen $\trop$ maps both the first and second quadrant to the region $R$ while it maps the third quadrant to the first quadrant.
So, it remains to show that an orbit cannot remain in the region $S$ indefinitely.
Consider $(s,t)$ in region $S$ so $\trop_1(s,t) = (-s, ps+t)$.
If $ps + t < 0$, then $\trop(s,t) = (-s, -ps - t)$ which is in the second quadrant.
If $ps + t \geq 0$, then $\trop(s,t) = (s', t') = ((pq-1)s + qt, -ps - t)$.
Assume that $(s',t')$ is in region $S$.
We have
\[\frac{t}{s} - \frac{t'}{s'} = \frac{ps^2 + pqst + qt^2}{ss'} > 0\]
and the change in angle is $\arctan\left( \frac{ps^2 + pqst + qt^2}{ss'+tt'} \right)$ by applying Lemma~\ref{lem:slope}.
If we stay in region $S$ and apply further iterations of $\trop$ this change in polar angle will only increase.
Indeed, $ps^2 + pqst + qt^2 = p(s')^2 + pqs't' + q(t')^2$ as this is the conserved quantity from Theorem~\ref{thm:conserved} while we also know that $s' < s$ and $|t'| < |t|$ given the direction of movement along the conserved quantity.
So, the change in angle is bounded away from zero.
Thus, we will eventually leave region $S$ implying we must reach region $R$.\\

If $(s,t)$ is in region $R$, then we have seen that $\trop(s,t)$ and hence $\trop^n(s,t)$ for all $n$ are in the region $R$.
Therefore we have established (ii).
Also, for such $(s,t)$ again letting $(s', t') = \trop(s,t)$, then we have
\begin{align*}
t' &= -ps - t\\
&\leq \sqrt{pq} t - t\\
&= (\sqrt{pq} - 1)t
\end{align*}
which is less than or equal to $t$ since $\sqrt{pq} \geq 2$.
This implies for $pq > 4$ the absolute value of the second coordinate will grow at least exponentially without bound when iterating $\trop$ and we have (iii) for $pq > 4$.
Now assume $pq = 4$, then $\sqrt{p} s + \sqrt{q} t > 0$ implies $ps + \sqrt{pq}t = ps + 2t > 0$.
Let $ps + 2t = \alpha > 0$ then
\begin{align*}
t' &= -ps -t\\
&= -(ps + 2t) + t\\
&= t - \alpha
\end{align*}
and the absolute value $|t'| = |t| + \alpha$ has grown.
Furthermore,
\begin{align*}
ps' + 2t' &= p((pq-1)s + qt) + 2(-ps - t)\\
&= p^2qs - ps + pqt - 2ps -2t\\
&= 4ps -ps + 4t - 2ps - 2t\\
&= ps + 2t\\
&= \alpha
\end{align*}
is conserved.
It follows that second coordinate will grow linearly without bound and (iii) is proven for $pq=4$.
The entire theorem is now proven.
\end{proof}

\begin{rmk}
When $p$ and $q$ are positive integers with $pq \geq 4$ statement (ii) in Theorem~\ref{thm:geq4} along with the analogous statement for $(\trop_1 \circ \trop_2)^n$ gives another proof of the rank 2 case of the asymptotic sign coherence conjecture~\cite[Theorem 3.5]{ASC}.
\end{rmk}

\bigskip

\section {Conjectures and Further Questions}

\bigskip

In this section we discuss some open problems on the dynamical system we have studied here.

\bigskip

We proved that the map $\mu$ from Section \ref{sec:x} is symplectic.
Additionally finding a conserved quantity would imply integrability of the map.
Although we have been unable to find one, the numerical
experiments suggest that a conserved quantity should exist. Therefore we would like to formulate the following open question.

\bigskip

\begin {problem} \label{prob:conserved}
    Find a conserved quantity for the map $\mu$.
\end {problem}

\bigskip

It was shown in~\cite{GNR} that each single mutation can be written as a finite-time evolution of a continuous Hamiltonian flow (see also~\cite{hkkr} for a similar phenomenon).
A simple computation shows that the same holds for $\mu_1$ and $\mu_2$ individually in our general situation.
\bigskip

\begin {problem}
Can the composition $\mu = \mu_2 \circ \mu_1$ be written as a finite-time evolution of a continuous Hamiltonian flow?
\end {problem}
\bigskip

We proved in Section \ref{subsec:unbounded} that the orbits are unbounded when $pq \geq 4$, but we did not show that the orbits are
bounded when $pq < 4$. 

\bigskip

\begin {problem} \label{prob:bounded}
    Show that the orbits of $\mu$ are bounded when $pq < 4$.
\end {problem}

\bigskip

\begin {rmk}
    It is possible that a solution to Problem \ref{prob:conserved}
    might imply \ref{prob:bounded} by studying the level sets of an explicit conserved quantity.
\end {rmk}

\bigskip

\section* {Acknowledgements}

\bigskip

N. Ovenhouse was supported by the NSF grant DMS-1745638. We would also like to thank Michael Shapiro for helpful comments on the content
and exposition of the paper.

\bibliographystyle{plainurl}
\bibliography{refs.bib}

\end {document}